\input amstex
\documentstyle{amsppt}
\magnification=\magstep1 \hoffset1 true pc \voffset2 true pc
\hsize36 true pc \vsize50 true pc \tolerance=2000

\define\gp#1{\langle#1\rangle}
\TagsOnRight \topmatter

\title
On units of  group algebras of $2$-groups of maximal class
\endtitle
\author Zs. Balogh, A. Bovdi
\endauthor

\leftheadtext\nofrills{Zs. Balogh, A. Bovdi}
\rightheadtext\nofrills{On units of  group algebras of $2$-groups
of maximal class}

\address
\hskip-\parindent Institute of Mathematics and Informatics,
College of Ny\'\i regyh\'aza,  S\'ost\'oi \'ut 31/b, H-4410 Ny\'\i
regyh\'aza, Hungary
\newline
baloghzs\@nyf.hu
\endaddress

\address
\hskip-\parindent Institute of Mathematics, University of
Debrecen, H-4010 Debrecen, P.O. Box 12, Hungary
\newline
bodibela\@math.klte.hu
\endaddress

\abstract We determine the number of elements of order two in the
group of normalized units $V({\Bbb F}_{2}G)$ of the group algebra
${\Bbb F}_{2}G$ of  a $2$-group of maximal class over the  field
${\Bbb F}_2$  of two elements. As a consequence for the $2$-groups
$G$ and $H$ of maximal class we have that $V({\Bbb F}_{2}G)$ and
$V({\Bbb F}_{2}H)$ are isomorphic if and only if $G$ and $H$ are
isomorphic.
\endabstract

\thanks
\newline Supported by Hungarian National Fund for Scientific
Research (OTKA) grants No.~T037202 and No.~T043034
\endthanks

\endtopmatter

\document
\subhead Notation and results
\endsubhead

Let $G$  be a finite  $p$-group  and ${\Bbb F}_{p}G$ its group
algebra over the field  ${\Bbb F}_{p}$ of $p$ elements. For
$x=\sum_{g \in G }\alpha_g g \in {\Bbb F}_{p}G $, \;$ \alpha_g\in
{\Bbb F}_{p}$ we  denote $\sum_{g \in G }\alpha_g$ by $\chi(x)$.
The subgroup
$$
V({\Bbb F}_{p}G)=\bigg\{\;  x=\sum_{g \in G }\alpha_g g \in {\Bbb
F}_{p}G  \quad\vert  \quad \chi(x)=1\; \bigg\}
$$
is called {\it the group of normalized units}. Evidently $V({\Bbb
F}_{p}G)$ is a $p$-group and its order is $p^{|G|-1}$.

\bigskip
 Let $C=\gp{\; a\;
\vert\; a^{2^n}=1\;}$ be  the  cyclic group of order $2^n$, where
$n\geq 2$. Consider the following extensions of $C$:

$$
\split Q_{2^{n+1}}&=\gp{\; a,b_1 \;\; \vert  \;\;
a^{2^n}=1, b_1^2=a^{2^{n-1}},  (a,b_1)=a^{-2},\; n\geq 2\;};\\
D_{2^{n+1}}&=\gp{\; a,b_2  \;\;\vert \;\;
a^{2^n}=1, b_2^2=1, (a,b_2)=a^{-2},\; n\geq 2\;};\\
SD_{2^{n+1}}&=\gp{\; a,b_3  \;\;\vert \;\; a^{2^n}=1, b_3^2=1,
(a,b_3)=a^{-2+2^{n-1}},\; n\geq 3\;},
\endsplit
$$
which are the   generalized quaternion group, the  dihedral group and
the semidihedral group respectively.
It is well-known that a finite  $2$-group of maximal class coincides with
one of these groups.

\newpage

\noindent Our main result is the following:

\proclaim {Theorem} Let $G$ be a $2$-group of maximal class and
let $\Theta_G(2)$ be the number of elements of order two in
$V(\Bbb F_2G)$. Then
$$\split \Theta_{D_{2^{n+1}}}(2)&=2^{2^n+n-1}+2^{2^n};\\
\Theta_{SD_{2^{n+1}}}(2)&=2^{2^n+n-1};\\
\Theta_{Q_{2^{n+1}}}(2)&=2^{2^n+n-1}-2^{2^n}. \endsplit$$
\endproclaim

The following question is due to S. D. Berman:  Let ${\Bbb F}_p$
be the field of $p$ elements, $G$ and $H$ finite $p$-groups. Is it
true that $V({\Bbb F}_{p}G)$ and $V({\Bbb F}_{p}H)$ are isomorphic
if and only if $G$ and $H$ are isomorphic?

It is an interesting and likely more difficult question than the
isomorphism problem  for modular group algebras of $p$-groups over
the field of $p$ elements.

In \cite{Berman} Berman gave a positive answer for this question
 for finite abelian $p$-groups. For
nonabelian  groups the question  is open. As a consequence of the
above theorem we have the

\proclaim {Corollary} Let ${\Bbb F}_2$ be the field of $2$
elements, and let $G$ and $H$ be finite $2$-groups of maximal
class. Then $V({\Bbb F}_{2}G)$ is isomorphic to $V({\Bbb F}_{2}H)$
if and only if $G$ and $H$ are isomorphic.
\endproclaim

Evidently if $V({\Bbb F}_{p}G)$ and  $V({\Bbb F}_{p}H)$ are not
isomorphic, then the group algebras ${\Bbb F}_{p}G$ and ${\Bbb
F}_{p}H$ are not isomorphic. Moreover it follows from Baginski's
result \cite{1} that the group algebra ${\Bbb F}_{2}G$  of a
$2$-group $G$ of maximal class uniquely determines the group $G$.

\bigskip

\subhead Some preliminary facts about group algebras of cyclic
groups
\endsubhead

\bigskip
We observe  some facts about  abelian $2$-groups  and their group
algebras.

Recall that for the  group algebra ${\Bbb F}_2C$ of a cyclic
$2$-group $C=\gp{\; a\, \vert\, a^{2^n}=1}$ the subset
$$
1,\,1+a,\,(1+a)^2,\dots,(1+a)^{2^n-1} \tag 1
$$
is a basis of ${\Bbb F}_2C$ and  a   normalized  unit $u$ can be written as
$$
1+\sum_{i=1}^{2^n-1}\alpha_i(a+1)^i.
$$
Each ideal of  the group algebra  ${\Bbb F}_2C$  has the form
${\Bbb F}_2C(1+a)^i$. Clearly $(1+a)^i$ belongs to the
augmentation ideal $A({\Bbb F}_2C)$ and by a Theorem of Jennings
\cite{6}
$$
A^i({\Bbb F}_2C)=A^{i+1}({\Bbb F}_2C)+(a+1)^i{\Bbb F}_2.
$$
According to \cite{5} the annihilator of the element $(1+a)^i$
coincides with $A^{2^n-i}({\Bbb F}_2C)={\Bbb F}_2C(1+a)^{2^n-i}$.

\proclaim{Lemma 1} The subgroup
$$
S_i=\{\; \gamma \in V({\Bbb F}_2C)  \quad \vert \quad
\gamma(1+a)^i=(1+a)^i\;\}$$ of $V({\Bbb F}_2C)$ has order $2^i$.
\endproclaim

\demo{Proof} It is easy to see that $1+Ann\big((1+a)^i\big)$
coincides with $S_i$, and according to  \cite{5} the cardinality
$|Ann((1+a)^i)|$ is $2^i$.

\rightline{$\square$}
\enddemo
\bigskip
Recall some well-known facts about a finite abelian $2$-group $N$.
Let $N[2]$ be  the subgroup generated by elements of order two,
and denote by  $\Phi(N)$  the Frattini subgroup of $N$. We shall
use the following

\proclaim{Lemma 2} Let $N$ be  a finite abelian $2$-group.
\item{$1$.} If $g \in \Phi(N)$, then $H_g=\{\;  h\in N\;\vert \;h^2=g
\;\}$ is a coset of $N$ by $N[2]$.
\item{$2$.} If $g \not\in \Phi(N)$, then $N=\gp{\; g}\times W$ for some
subgroup $W$.
\endproclaim

\bigskip\bigskip

\subhead Involutions and unitary subgroups
\endsubhead

\bigskip
Let $C=\gp{\; a  \;\vert \; a^{2^n}=1,\; n\geq 2\;}$ be a cyclic
$2$-group. First we review some results on involutions of ${\Bbb
F}_{2}C$. Recall that the linear extension of the automorphism
$a^i \rightarrow a^{-i}$ of $G$ to the automorphism $x\rightarrow
x^*$ of ${\Bbb F}_{2}C$ is an involution of ${\Bbb F}_{2}C$.
Moreover, for $n\geq 3$ we have another involution $x\rightarrow
x^{\circledast}$, which is generated by the automorphism $a^i
\rightarrow a^{(2^{n-1}-1)i}$ of order two. For the sake of
convenience we assume that $\sigma$ is either $x\rightarrow x^*$
or $x\rightarrow x^{\circledast}$.

First we observe some properties of these involutions.

\proclaim{Lemma 3} For
$x=\alpha_0+\alpha_1a+\dots+\alpha_{2^n-1}a^{2^n-1} \in {\Bbb
F}_{2}C$ we have
$$
x^2=\sum_{i=0}^{2^{n-1}-1}(\alpha_i+\alpha_{i+2^{n-1}})a^{2i}
$$
and
$$
xx^*=\chi(x)+\sum_{j=1}^{2^{n-1}-1}\Bigl[ \sum_{i=0}^{2^n-1}
\alpha_i\alpha_{i-j \pmod{2^n}} \Bigr]\cdot\Bigl[a^j+a^{-j}\Bigr].
$$
\endproclaim

\demo{Proof} Evidently $x^2=\sum_{i=0}^{2^n-1}\alpha_ia^{2i}$ and
$(a^{i+2^{n-1}})^2=(a^{i})^2$, thus the coefficient of $a^{2i}$
equals $\alpha_i+\alpha_{i+2^{n-1}}$. Moreover, the equality
$(xx^*)^*=xx^*$ asserts that
$$
xx^*=\gamma_0+\gamma_{2^{n-1}}a^{2^{n-1}}+\sum_{j=1}^{2^{n-1}-1}
\gamma_j(a^j+a^{-j}).
$$
For each $a^j \in supp\,(xx^*)$ with $0\leq j\leq 2^{n-1}-1$ the
equality $a^j=a^ia^{-(i-j)}$ shows that
$\gamma_j=\sum_{i=0}^{2^n-1} \alpha_i\alpha_{i-j \pmod{2^n}}$.

In particular, $\gamma_0=\sum_{i=0}^{2^n-1} \alpha_i\alpha_i
=\sum_{i=0}^{2^n-1} \alpha_i= \chi(x)$ and
$$
\gamma_{2^{n-1}}= \sum_{i=0}^{2^n-1}\alpha_i\alpha_{i-2^{n-1}
\pmod{2^n}}= 2\cdot
\sum_{i=0}^{{2^{n-1}}-1}\alpha_i\alpha_{i-2^{n-1}
\pmod{2^n}}=0.\quad \quad\square
$$
\enddemo
\noindent We put
$$
Q=\{\; 0,2,4,\ldots,2^{n-2}-2\;\} \cup \{\;
2^{n-1},2^{n-1}+2,2^{n-1}+4,\ldots,2^{n-1}+2^{n-2}-2 \;\},
$$
$P=\{\; 0,2,4,\ldots,2^{n}-2\;\}$  and $R=\{\;
0,2,4,\ldots,2^{n-1}-2\;\}$.

\bigskip

\proclaim{Lemma 4} For
$x=\alpha_0+\alpha_1a+\dots+\alpha_{2^n-1}a^{2^n-1} \in {\Bbb
F}_{2}C$ we have
$$
xx^{\circledast}
  =\gamma_0+\gamma_{{2^{n-1}}}a^{2^{n-1}}+\sum_{k\in R\setminus \{\; 0\;\}}\gamma_{k}(a^k+a^{-k})
  +\sum_{k\in 1+Q}\gamma_{k}(a^k+a^{-k+2^{n-1}}),
$$
where
$$
\gamma_{k}=\cases\sum\limits_{r \in P}\alpha_r\alpha_{r-k
\pmod{2^n}}+ \sum\limits_{r \in 1+P}\alpha_r
\alpha_{r-k+2^{n-1} \pmod{2^n}} \;  \text{for }    k \in P\setminus\{\; 0,2^{n-1}\;\};\\
\sum\limits_{r \in P} \alpha_r \alpha_{r-k+2^{n-1} \pmod{2^n}} +
\sum\limits_{r \in 1+P} \alpha_r\alpha_{r-k \pmod{2^n}} \;
\text{for }    k \in 1+P;\\
\sum\limits_{r \in 1+P}\alpha_r \quad \text{for  }
k=2^{n-1};  \\
\sum\limits_{r \in P}\alpha_r \quad \quad \text{for  } k=0.
\endcases
$$
\endproclaim

\demo{Proof} Since $(xx^{\circledast})$ is $\circledast$-symmetric
so
$$
xx^{\circledast}
  =\gamma_0+\gamma_{{2^{n-1}}}a^{2^{n-1}}+\sum_{k\in R\setminus\{\; 0\;\}}\gamma_{k}(a^k+a^{-k})
  +\sum_{k\in 1+Q}\gamma_{k}(a^k+a^{-k+2^{n-1}}).
$$
We define the permutation $\rho$ of the set $\{\;
0,1,2,\dots,2^n-1 \;\}$ in the following way:
$$
\rho(i)=\cases i& \text{ if $i$ is even};\\
i+2^{n-1} \pmod{2^n}& \text{ if $i$ is odd}.
\endcases
$$
Using the permutation $\rho$, simple computations show that
$x^{\circledast}=\sum_{i=0}^{2^{n}-1} \alpha_{\rho(i)} a^{-i}$ and
the coefficient of $a^k \in supp\,(xx^{\circledast})$ is equal to
the trace $\gamma_k=tr(xx^{\circledast}a^{-k})$. Evidently
$$
\gamma_k=tr\Big(\big(\sum_{j=0}^{2^{n}-1} \alpha_j a^j\big)\cdot
\big(\sum_{i=0}^{2^{n}-1} \alpha_{\rho(i-k)} a^{-i}\big)\Big)
=\sum_{r=0}^{2^{n}-1}\alpha_r\alpha_{\rho(r-k)}.
$$
Therefore we have for even $k$
$$
\gamma_{k}=\sum_{r=0}^{2^{n}-1}\alpha_r\alpha_{\rho(r-k)}=\sum\limits_{r
\in P}\alpha_r\alpha_{r-k \pmod{2^n}} +\sum\limits_{r \in 1+P}
\alpha_r\alpha_{r-k+2^{n-1} \pmod{2^n}},
$$
and for odd $k$
$$
\gamma_{k}=\sum\limits_{i \in P} \alpha_i \alpha_{i-k+2^{n-1}
\pmod{2^n}} + \sum\limits_{i \in 1+P} \alpha_i\alpha_{i-k
\pmod{2^n}}.
$$

\noindent Specifically for $k=2^{n-1}$ simple computations show
that
$$
\gamma_{{2^{n-1}}}=\sum\limits_{i \in P }\alpha_i\alpha_{i-2^{n-1}
\pmod{2^n}}+
  \sum\limits_{i \in 1+P}\alpha_i\alpha_{i}=\sum\limits_{i \in 1+P}\alpha_i,
$$
and
$$
\gamma_0=\sum\limits_{i \in P}\alpha_i\alpha_i+ \sum\limits_{i \in
1+P}\alpha_i\alpha_{i-2^{n-1} \pmod{2^n}}= \sum\limits_{i \in
P}\alpha_i.
$$

\rightline{$\square$}
\enddemo

\bigskip
Recall that each involution $\sigma$ of ${\Bbb F}_2C$ determines a
$\sigma$-unitary subgroup
$$
V_{\sigma}({\Bbb F}_2C)=\{\; \gamma \in V({\Bbb F}_2C)  \quad
\vert \quad \gamma^{-1}=\gamma^{\sigma} \;\}
$$
of $V({\Bbb F}_2C)$.

The structure of $*$-unitary subgroup of $V({\Bbb F}_2C)$ was
described in \cite{3} and \cite{4}. According to these results the
order of $V_*({\Bbb F}_2C)$ is
$$
|C^2[2]|\cdot|{\Bbb
F}_2|^{\frac{1}{2}(|C|+|C[2]|)-1}=2^{2^{n-1}+1}. \tag 2
$$
Now we determine the order of the $\circledast$-unitary subgroup
$V_{\circledast}({\Bbb F}_2C)$, and as far as we know, it has not
been investigated so far.

The mapping $\varphi_{\sigma}$, given by
$$
\varphi_{\sigma}(x)=x^{\sigma}x^{-1} \quad \text{for }  \quad x
\in V({\Bbb F}_2C),
$$
is a homomorphism of the group $V({\Bbb F}_2C)$ onto some subgroup
$W_{\sigma}(C)$ of the \newline $\sigma$-unitary subgroup
$V_{\sigma}({\Bbb F}_2C)$ with kernel $S_{\sigma}(C)=\{\;  x\in
V({\Bbb F}_2C) \quad \vert  \quad x=x^{\sigma} \;\}$.

The subset $D$ of $C$  determines the element $\widehat{D}=\sum_{g
\in D} g$ of $F_2C$.

\proclaim{Lemma 5} The unit $1+\widehat C$ does not belong to
$W_{\circledast}(C)$ and $1+\widehat{C^2}$ is not an element of the
subgroup $V_{\circledast}^2({\Bbb F}_2C)$.
\endproclaim

\demo{Proof} Assume that $1+\widehat C \in W_{\circledast}(C)$.
Then $1+\widehat C=x^{\circledast}x^{-1}$ for some $x \in V({\Bbb
F}_2C)$ and
$$
x^{\circledast}=x(1+\widehat C)=x+\chi(x)\widehat C=x+\widehat
C.\tag 3
$$
Evidently the traces of the elements $x$ and $x^{\circledast}$ are
equal, so $tr(x+x^{\circledast})=0$ and the equality $(3)$ leads
to contradiction at characteristic two.

Now, suppose that $1+\widehat{C^2}=x^2$ for some
$\circledast$-unitary element $x \in V_{\circledast}({\Bbb
F}_2C)$. Since $x^{\circledast}=x^{-1}$ and $x$ has the form
$x=y_1+y_2a$ with $y_i \in {\Bbb F}_2C^2$, so
$$
x^{\circledast}+x=x(1+x^{-2})=x(1+x^2)= \widehat{C^2}\Big(\chi(y_1)+\chi(y_2)a^{\circledast}\Big),
$$
and $\chi(y_1)=tr(x^{\circledast}+x)=0$. This shows that
$y_1^{\circledast}+y_1+y_2a+(y_2a)^{\circledast}=\widehat{C^2}a^{\circledast}$,
which implies that $y_1=y_1^{\circledast}$ and
$$
\widehat{C^2}=y_2^{\circledast}+y_2a^{2^{n-1}+2}.
$$
The previous equality shows that the set
$supp\;(y_2^{\circledast}) \cap supp\;(y_2a^{2^{n-1}+2})$ is
empty, so $|supp\;(y_2^{\circledast})|=2^{n-2}$, thus
$\chi(y_2)=0$, a contradiction.

\rightline{$\square$}
\enddemo

The next homomorphism of $V({\Bbb F}_2C)$ we shall use later.
Evidently if $\Psi_{\sigma}(x)=xx^{\sigma}$ for $x \in V({\Bbb
F}_2C)$, then $\Psi_{\sigma}$ is a homomorphism of $V({\Bbb
F}_2C)$ with kernel $V_{\sigma}({\Bbb F}_2C)$.

\proclaim{Lemma 6} If $C=\gp{\; a \;\vert \; a^{2^n}=1,\; n\geq
3\; }$, then the order of $V_{\circledast}({\Bbb F}_{2}C)$ is
$2^{\frac{|C|}{2}}$ and
$$
V_{\circledast}({\Bbb F}_{2}C)=W_{\circledast}(C) \times \gp{\;
1+\widehat C}.
$$
\endproclaim

\demo{Proof} It is easy to see that the lower layers of the both
groups $V_{\circledast}({\Bbb F}_2C)$ and $S_{\circledast}(C)$ coincide.
 To determine the $2$-rank of $V_{\circledast}({\Bbb
F}_2C)$ it suffices to find the order of the lower layer
$S_{\circledast}(C)[2]$ of the group $S_{\circledast}(C)$.

Each $\circledast$-symmetric element of $V({\Bbb F}_2C)$ has the
form
$$
\gamma_0+\gamma_{{2^{n-1}}}a^{2^{n-1}}+\sum_{k\in R\setminus\{\;
0\;\}}\gamma_{k}(a^k+a^{-k}) +\sum_{k\in
1+Q}\gamma_{k}(a^k+a^{-k+2^{n-1}}),
$$
where $\gamma_{i}\in {\Bbb F}_2$ and
$\gamma_0+\gamma_{{2^{n-1}}}=1$. Therefore
$|S_{\circledast}(C)|=2^{\frac{|C|}{2}}$ and analogously,
$|S_{\circledast}(C)^2|=2^{\frac{|C^2|}{2}-1}$. This asserts that
the order of the subgroups $S_{\circledast}(C)[2]$ and
$V_{\circledast}({\Bbb F}_2C)[2]$ is $2^{\frac{|C|}{4}+1}$. But
the kernel of $\varphi_{\circledast}$ is $S_{\circledast}(C)$,  so
$$
|W_{\circledast}(C)|=|V({\Bbb
F}_2C)|:|S_{\circledast}(C)|=2^{\frac{|C|}{2}-1}.
$$

Now  we are ready to prove by induction on the order of $C$ that
$|V_{\circledast}({\Bbb F}_{2}C)|=2^{\frac{|C|}{2}}$. First let
$|C|=2^3$. According to Lemma $4$ for  the unit
$x=\sum_{i=0}^{7}\alpha_{i} a^i$ we have
$$
xx^{\circledast}=\beta_0+\beta_1(a+a^3+a^5+a^7)+\beta_2(a^2+a^6)+
(\beta_{0}+1)a^{4},
$$
because $\chi(x)=1$.

Since $\psi_{\circledast}(x)=xx^{\circledast}$, and the order of
$Im(\psi_{\circledast})$ is $2^3$, so the order of
$V_{\circledast}({\Bbb F}_2C)$ is
$$
|V({\Bbb F}_2C)|:|Im(\psi_{\circledast})|=2^4.
$$
Evidently $1+\widehat C \in V_{\circledast}({\Bbb F}_2C)$ and
$1+\widehat C \not\in W_{\circledast}(C)$ by Lemma $5$. Since the
subgroup $\gp{\; 1+\widehat C} \times W_{\circledast}(C)$ of
$V_{\circledast}({\Bbb F}_2C)$ has order $2^4$ thus
$$
V_{\circledast}({\Bbb F}_2C)=\gp{\; 1+\widehat C} \times
W_{\circledast}(C).
$$
Now let $|C|>2^3$. Applying the inductive hypothesis for $C^2$,
$$
V_{\circledast}({\Bbb F}_2C^2)=W(C^2)\times\gp{\;
1+\widehat{C^2}}.
$$
Lemma $5$ asserts that $\gp{\; 1+\widehat C} \times
W_{\circledast}(C) $ is a subgroup of $V_{\circledast}({\Bbb
F}_2C)$ and \newline $\gp{\; 1+\widehat{C^2}} \times
V_{\circledast}^2({\Bbb F}_2C)$ is a subgroup of
$V_{\circledast}({\Bbb F}_2C^2)$. This shows that
$$
|V_{\circledast}^2({\Bbb F}_2C)|\leq | V_{\circledast}({\Bbb
F}_2C^2)|:|\gp{\; 1+\widehat{C^2}}|=|W(C^2)|,
$$
so
$$
|V_{\circledast}({\Bbb F}_2C)|=|V_{\circledast}^2({\Bbb
F}_2C)|\cdot |V_{\circledast}({\Bbb F}_2C)[2]|\leq |W(C^2)|\cdot
|V_{\circledast}({\Bbb F}_2C)[2]|=2^{\frac{|C|}{2}}.
$$
But $\gp{\; 1+\widehat C} \times W_{\circledast}(C)$ is a subgroup
of $V_{\circledast}({\Bbb F}_2C)$ and its order is
$2^{\frac{|C|}{2}}$, therefore the order of $V_{\circledast}({\Bbb
F}_2C)=\gp{\; 1+\widehat C} \times W_{\circledast}(C)$ is
$2^{\frac{|C|}{2}}$.

\rightline{$\square$}
\enddemo
\bigskip

\subhead Elements of order two in $V({\Bbb
F}_{2}G)$
\endsubhead
\bigskip

For a fixed noninvertible element $z \in {\Bbb F}_{2}C$ the set
$$
M_{z}^{\sigma}=\{\; y\in V({\Bbb F}_{2}C) \,|\, (y+y^{\sigma})z=0
\;\}
$$
is a subgroup of $V({\Bbb F}_{2}C)$. Indeed, if $y_1,y_2 \in
M_{z}^{\sigma}$, then $y_1z=y_1^{\sigma}z$ and
$y_2z=y_2^{\sigma}z$. Hence
$y_1y_2z=y_1y_2^{\sigma}z=y_1^{\sigma}y_2^{\sigma}z=(y_1y_2)^{\sigma}z$,
so $y_1y_2 \in M_{z}^{\sigma}$.

In $V({\Bbb F}_{2}G)$ we divide the units of order two into two
classes. It is well-known that $x=x_1+x_2b$ is a unit if and only if
$\chi(x_1)+\chi(x_2)=1$.

\proclaim{Definition} Let $x=x_1+x_2b$ be a unit. If $x$ has order
two, $\chi(x_1)=1$ and $\chi(x_2)=0$, then $x$ is called {\it a
unit of type $1$}. If $x$ has order two, $\chi(x_1)=0$ and
$\chi(x_2)=1$, then $x$ is called {\it a unit of type $2$}.
\endproclaim

\proclaim{Lemma 7} The number of units of type $1$ in both groups
$V(F_2D_{2^{n+1}})$ and $V(F_2Q_{2^{n+1}})$ is equal.
\endproclaim

\demo{Proof}  A unit $x_1+x_2b_2 \in V({\Bbb F}_{2}D_{2^{n+1}})$
has order two if and only if
$$\cases
x_1^2=x_2x_2^*+1;\\
(x_1+x_1^*)x_2=0.\endcases\tag 4
$$
Similarly, $x_1+x_2b_1 \in V({\Bbb F}_{2}Q_{2^{n+1}})$ is a unit
 of order two if and only if
$$\cases
x_1^2=x_2x_2^*a^{2^{n-1}}+1;\\
(x_1+x_1^*)x_2=0.
\endcases
$$
\vskip8pt

Let $x_2 \in {\Bbb F}_2C$ be a fixed not invertible element.
Clearly $x_2x_2^*+1$, $x_2x_2^*a^{2^{n-1}}+1$ belong to the
subgroup $M_{x_2}$, defined before, and the set
$$
H_{x_2x_2^*+1}=\{\; \; h \in M_{x_2}^*\quad \vert\quad
h^2=x_2x_2^*+1\;\;\}
$$
either is empty or according to Lemma $2$ constitutes a coset of
$M_{x_2}^*$ by $M_{x_2}^*[2]$. Similarly,
$$
H_{x_2x_2^*a^{2^{n-1}}+1}=\{\; \;h \in M_{x_2}^* \quad \vert\quad
h^2=x_2x_2^*a^{2^{n-1}}+1\;\;\}
$$
is either a coset of $M_{x_2}^*$ by $M_{x_2}^*[2]$ or an empty
set.

Let us prove that $x_2x_2^*+1 \in \Phi(M_{x_2}^*)$ if and only if
$\quad x_2x_2^*a^{2^{n-1}}+1 \in \Phi(M_{x_2}^*)$.

\noindent Suppose that $u=x_2x_2^*+1 \notin \Phi(M_{x_2}^*)$ and
$v=x_2x_2^*a^{2^{n-1}}+1 \in \Phi(M_{x_2}^*)$. Then for $u \notin
\Phi(M_{x_2}^*)$ there exists a direct decomposition
$M_{x_2}^*=\gp{\; u} \times W$ for some subgroup $W$. Of course
$$
v^2=(x_2x_2^*a^{2^{n-1}}+1)^2=(x_2x_2^*)^2+1=u^2,
$$
therefore $u^{-1}v\in M_{x_2}^*[2]$ and $v=uy$ for some $y\in
M_{x_2}^*[2]$. Obviously,  $v \in \Phi(M_{x_2}^*)=\gp{\;
u^2}\times W^2$, whence $v=u^{2t}w^2$ for some $t$ and $w\in W$.
Moreover, $v=uy$ implies that
$$
y=u^{2t-1}w^2 \in M_{x_2}^*[2]\quad \text{and}\quad
1=y^2=u^{4t-2}w^4.
$$
This shows that $u^{4t-2}=1$,  so $4t \equiv
2 \pmod{|u|}$, where $|u|$ is the order of $u$. Hence $|u|=2$ and
$(x_2x_2^*)^2=0$. By Lemma $3$, $x_2x_2^*=z(1+a^{2^{n-1}})$ for
some $z \in {\Bbb F}_{2}C$ and
$$
v=x_2x_2^*a^{2^{n-1}}+1=z(1+a^{2^{n-1}})a^{2^{n-1}}+1=z(1+a^{2^{n-1}})+1=x_2x_2^*+1=u
$$
which is impossible.

We have shown that the cardinalities of the subsets
$H_{x_2x_2^*+1}$ and $H_{x_2x_2^*a^{2^{n-1}}+1}$ are equal for
each not invertible element $x_2$, thus the proof of lemma is
complete.

\rightline{$\square$}
\enddemo

A  unit $x_1+x_2b_3 \in V({\Bbb F}_{2}SD_{2^{n+1}})$ has order
two if and only if
$$\cases
x_1^2=x_2x_2^{\circledast}+1;\\
(x_1+x_1^{\circledast})x_2=0.\endcases\tag 5
$$
\noindent For each noninvertible and nonzero element $x_2$ there
exists an $i$ such that \newline $x_2=\gamma(1+a)^i$ for some
$\gamma \in V({\Bbb F}_2C)$. The equalities $(4)$ and $(5)$ assert
that
$$
x_2x_2^*=\gamma\gamma^*(1+a)^{2i}a^{-i}\in {\Bbb F}_2C^2,
$$ and
$$
x_2x_2^{\circledast}=\gamma\gamma^{\circledast}(1+a)^{i}(1+a^{2^{n-1}-1})^i\in
{\Bbb F}_2C^2.
$$
In particular, if $i=2l$ is even, then
$$
(1+a)^{2l}(1+a^{2^{n-1}-1})^{2l}=(1+a^2)^{l}(1+a^{-2})^{l}=(1+a)^{2l}(1+a^{-1})^{2l},
$$
and \quad
$x_2x_2^{\circledast}=\gamma\gamma^{\circledast}(1+a)^{2i}a^{-i}$.

For each  $0 \leq i< 2^n$ we define the set
$$
H_i^{\sigma}=\{\;  \; h \in V({\Bbb F}_2C) \quad \vert\quad
hh^{\sigma}(1+a)^i(1+a^{\sigma})^i\in {\Bbb F}_2C^2 \; \;\}.
$$

\proclaim{Lemma 8} The set $H_i^{\sigma}$ has the following
properties:
\item{$1$.} If $i\geq 2^{n-1}$, then $H_i^{\sigma}=V({\Bbb F}_2C)$.
\item{$2$.} If $i=2l+1$ is odd and $i<2^{n-1}$, then the set
$H_i^{\sigma}$is empty.
\item{$3$.} If $i=2l$ is even and $i<2^{n-1}$, then
 $H_{2l}^{\sigma}$ is a subgroup of $V({\Bbb F}_2C)$.
\item{$4$.} For even indices $i=2l<2^{n-1}$
$$
H_0^{\sigma}\subset H_2^{\sigma}\subset H_4^{\sigma} \subset \cdots
\subset H_{2l}^{\sigma} \subset \cdots \subset
H_{2^{n-1}-2}^{\sigma}=V({\Bbb F}_2C)
$$
and the order of the subgroup $H_{2l}^{\sigma} $ is $2^{3\cdot
2^{n-2}+l}$.
\endproclaim

\demo{Proof} {\bf $1$.} Since
$(1+a)^{2^{n-1}}=(1+a^{\sigma})^{2^{n-1}}$, for $i\geq 2^{n-1}$ we
have
$$
hh^{\sigma}(1+a)^i(1+a^{\sigma})^i
=hh^{\sigma}(1+a)^{2^{n}}(1+a)^{i-2^{n-1}}(1+a^{\sigma})^{i-2^{n-1}}=0.
$$
Thus  $hh^{\sigma}(1+a)^i(1+a^{\sigma})^i\in {\Bbb F}_2C^2$ for
each  $h\in V({\Bbb F}_2C)$.
Hence $H_i^{\sigma}=V({\Bbb F}_2C)$ for all $i\geq 2^{n-1}$.

{\bf $2$.} Let $i=2l+1$ be a fixed odd integer and $i<2^{n-1}$.

\noindent First, let us prove by induction on $l$ that if
$i=2l+1$, then
$$
(a+a^{-1})^i=\sum_{\text{r is odd}}\beta_r (a^r+a^{-r}),  \tag 6
$$
where $\beta_r \in {\Bbb F}_2$. It is clear for $l=0$ and assume
that
$$
(a+a^{-1})^{2l+1}=\sum_{\text{r is odd}}\beta_r (a^r+a^{-r}).
$$
By the identity
$$
(a^s+a^{-s})(a^k+a^{-k})=(a^{s+k}+a^{-(s+k)})+(a^{s-k}+a^{-(s-k)})
\tag 7
$$
we have
$$
(a+a^{-1})^{2(l+1)+1}=\sum_{\text{r is odd}}\beta_r
\big((a^{r+2}+a^{-r-2})+(a^{r-2}+a^{-r+2})\big)
$$
the desired assertion.

Now, we begin to prove that the set $H_i^{\sigma}$ is empty. First
let $\sigma$ be the $*$-involution. By Lemma $3$, for
$h\in H_i^{*}$ we have  $ hh^*=1+z_1+z_2$, where
$$
z_1=\sum_{j \in R\setminus\{\; 0\;\}}\alpha_{j}(a^j+a^{-j}), \quad
z_2=\sum_{k \in 1+R}\alpha_{k}(a^k+a^{-k}),
$$
and $R=\{\; 0,2,4,\ldots, 2^{n-1}-2\;\}$. By definition,
$$
hh^*(a+a^{-1})^i=(1+z_1)(a+a^{-1})^i+z_2(a+a^{-1})^i \in {\Bbb
F}_2C^2\tag{8}
$$
and since $i=2l+1$ is odd it is easy to see that
$$
(a^k+a^{-k})(a+a^{-1})^i=(a^{k}+a^{-k})(a+a^{-1})(a^2+a^{-2})^l,
$$
so by $(6)$ and $(7)$ we have $z_2(a+a^{-1})^i\in {\Bbb F}_2C^2$
and $w=(1+z_1)(a+a^{-1})^i \not \in {\Bbb F}_2C^2$ if $w\ne 0$.
But $(8)$ confirms that $(1+z_1)(a+a^{-1})^i=0$ and we have
$(a+a^{-1})^i=0$, because $1+z_1$ is a unit. This is impossible
for $i<2^{n-1}$.

Let $\sigma$ be the $\circledast$-involution. Lemma $4$ asserts
that $hh^{\circledast}=1+z_1+z_2$, where
$$
z_1=(\gamma_0+1)(1+a^{2^{n-1}})+\sum_{\sixrm k\in R\setminus\{\;
0\;\}}\gamma_{k}(a^k+a^{-k}),
$$
$$
z_2=\sum_{k\in 1+Q}\gamma_{k}(a^k+a^{-k+2^{n-1}})
$$
and
$$
Q=\{\;  0,2,4,\ldots,2^{n-2}-2 \;\}\cup \{\;  2^{n-1},2^{n-1}+2
,2^{n-1}+4,\ldots 2^{n-1}+2^{n-2}-2\;\}.
$$
We remark that
$$
hh^{\circledast}(1+a)^{2l+1}(1+a^{2^{n-1}-1})^{2l+1} =
$$$$
=(1+z_1+z_2)(a+a^{-1})^{2l}\Big((1+a)^{2^{n-1}}+(a+a^{-1})+(1+a)^{2^{n-1}}a^{-1}\Big)
\in {\Bbb F}_2C^2.
$$
The identity
$$
(a^k+a^{-k+2^{n-1}})(a^r+a^{-r})=(a^{k+r}+a^{-r-k+2^{n-1}})+(a^{k-r}+a^{-k+r+2^{n-1}})
$$
and $(7)$ assert that
$$
(1+z_1)(a+a^{-1})^{2l}(1+a)^{2^{n-1}}+z_2(a+a^{-1})^{2l}(a+a^{-1})+
z_2(a+a^{-1})^{2l}(1+a)^{2^{n-1}}a^{-1}
$$
belongs to ${\Bbb F}_2C^2 $, and
$$
w=(1+z_1)(a+a^{-1})^{2l+1}+(a+a^{-1})^{2l}(1+a)^{2^{n-1}}z_2+
(a+a^{-1})^{2l}(1+a)^{2^{n-1}}(1+z_1)a^{-1}
$$
is not an element of ${\Bbb F}_2C^2$, if  $w\neq 0$. By $(8)$, $w \in
{\Bbb F}_2C^2$, but this leads again to a contradiction. Therefore
$$
(1+z_1)(a+a^{-1})^{2l+1}+(a+a^{-1})^{2l}(1+a)^{2^{n-1}} \Big( z_2+
(1+z_1)a^{-1}\Big)=0.
$$

\noindent Since $1+z_1$ and $e=z_2(1+z_1)^{-1}+a^{-1}$ are units,
it follows that
$$
(a+a^{-1})^{2l}\Big((a+a^{-1})+(1+a)^{2^{n-1}}e\Big)=0.\tag 9
$$

For $l=0$ we have $(a+a^{-1})=(1+a)^{2^{n-1}}e$, which is
impossible, because $(a+a^{-1})^2\ne 0$. If $l>0$, then
$$
(a+a^{-1})^{2l}(1+a)^{2^{n-1}}e \in A^{4l+2^{n-1}}({\Bbb F}_2C).
$$
Thus from $(9)$ it follows that $(a+a^{-1})^{i} \in
A^{4l+2^{n-1}}({\Bbb F}_2C)$ and $i<4l+2^{n-1}$. This is
impossible by Jenning's theory \cite{4}. Thus $H_i^{\sigma}$ is
empty, as we stated.
\bigskip

{\bf $3$.} Now we shall prove that for $i=2l<2^{n-1}$ the set
$H_i^{\sigma}$ is a subgroup of $V({\Bbb F}_2C)$. For $h \in H_i^{\sigma}$
using the basis
$(1)$ we  have
$$
hh^{\sigma}=1+\sum_{j=1}^{2^n-1}\alpha_j(1+a)^j,
$$
and
$$
hh^{\sigma}(1+a)^{2i}=(1+a)^{2i}+\sum_{j=1}^{2^n-2i-1}\alpha_j(1+a)^{j+2i}
\in {\Bbb F}_2C^2.
$$
But it follows that
$\alpha_1=\alpha_3=\alpha_5=\dots=\alpha_{2^n-2i-1}=0$.
Consequently, for $h_1,h_2 \in H_i^{\sigma}$ there exist $u_1,u_2
\in {\Bbb F}_2C^2$ and $v_1,v_2 \in A^{2^n-2i}({\Bbb F}_2C)$ such
that
$$
h_1h_1^{\sigma}= 1+u_1+v_1 \quad \text{and} \quad h_2h_2^{\sigma}
= 1+u_2+v_2.
$$
Clearly
$$
h_1h_2(h_1h_2)^{\sigma}=h_1h_1^{\sigma}h_2h_2^{\sigma}=1+u_1+u_2+u_1u_2+z,
$$
where $1+u_1+u_2+u_1u_2 \in {\Bbb F}_2C^2$ and
$z=v_1+v_1u_2+v_2+u_1v_2+v_1v_2 \in A^{2^n-2i}({\Bbb F}_2C)$.
Hence $h_1h_2 \in H_i^{\sigma}$ and $H_i^{\sigma}$ is a subgroup
of $V({\Bbb F}_2C)$.
\bigskip

In the next step we shall verify that
$H_{2^{n-1}-2}^{\sigma}=V({\Bbb F}_2C)$. A  unit $u$ of ${\Bbb
F}_2C$ can be written as
$$
u\equiv 1+\alpha_1(1+a)+\alpha_2(1+a)^2+\alpha_3(1+a)^3
\pmod{A^4({\Bbb F}_2C)}.
$$
It is easy to see that
$$
\split (1+a^{\sigma}) &\equiv(1+a)+(1+a)^2+\\
&+(1+a)^3+\dots+(1+a)^{2^{n-1}-1} \pmod{A^{2^{n-1}}({\Bbb
F}_2C)}.\\
\endsplit \tag 10
$$
It follows that
$$
u^{\sigma}\equiv
1+\alpha_1(1+a)+(\alpha_1+\alpha_2)(1+a)^2+(\alpha_1+\alpha_3)(1+a)^3
\pmod{A^4({\Bbb F}_2C)}
$$
and $$uu^{\sigma} \equiv 1 \pmod{A^4({\Bbb F}_2C)}. \tag 11$$

\noindent Therefore for each $h \in V({\Bbb F}_2C)$ and
$i=2^{n-1}-2$, $(11)$ confirms that
$$
hh^{\sigma}(1+a)^{2i}=hh^{\sigma}(1+a)^{(2^n-4)}=(1+a)^{(2^n-4)}.
$$
This shows that $h\in H_{2^{n-1}-2}^{\sigma}$ and
$H_{2^{n-1}-2}^{\sigma}=V({\Bbb F}_2C)$, because $h$ is an
arbitrary unit. This completes the proof of this assertion.
\bigskip

{\bf $4$.} Now let $i=2l< 2^{n-1}-2$. We shall prove that
$H_{i}^{\sigma}$ is a proper subgroup of $H_{i+2}^{\sigma}$.
Clearly $H_{i}^{\sigma} \subseteq H_{i+2}^{\sigma}$ and it is
sufficient to verify that $h=1+(1+a)^{2^n-(2i+3)}$ does not belong
to $H_i^{\sigma}$ and $h \in H_{i+2}^{\sigma}$.

\noindent Note that for an even $i$ the binomial coefficient
$\textstyle\big({2^n-2i-3 \atop 2}\big)$ is even too and
\newline $h^{\sigma}=1+(1+a^{\sigma})^{2^n-(2i+3)}$. By the
binomial formula and $(10)$ we have
$$
h^{\sigma}=
1+(1+a)^{2^n-(2i+3)}+(1+a)^{2^n-(2i+2)}+(1+a)^{2^n-(2i+1)}+x_{\sigma},
$$
where $x_{\sigma} \in A^{2^n-2i}({\Bbb F}_2C)$. It is easy to see that
$2^{n+1}-4i-j\geq 2^n-2i$ for  $i\leq
2^{n-1}-4$ and
$j=4,5,6$. This shows that $ hh^{\sigma}$ can be written as
$$
hh^{\sigma}=1+(1+a)^{2^n-(2i+2)}+(1+a)^{2^n-(2i+1)}+y_{\sigma}
$$
for some $y_{\sigma} \in A^{2^n-2i}({\Bbb F}_2C)$. Thus
$$
hh^{\sigma}(1+a)^{2i}=(1+a)^{2i}+(1+a)^{2^n-2}+(1+a)^{2^n-1}
$$
and we have  $h \notin H_i^{\sigma}$. But the equality
$hh^{\sigma}(1+a)^{2i+2}=(1+a)^{2i+2}$ states that $h \in
H_{i+2}^{\sigma}$ and this proves that
$$
H_0^{\sigma}\subset H_2^{\sigma}\subset H_4^{\sigma} \subset \dots
\subset H_i^{\sigma} \subset \dots \subset
H_{2^{n-1}-2}^{\sigma}=V({\Bbb F}_2C).
$$
\smallskip

Now we shall determine the order of the group $H_0^{\sigma}$. The
$\sigma$-unitary subgroup $V_{\sigma}({\Bbb F}_2C)$ of $V({\Bbb
F}_2C)$ is a subgroup of $H_0^{\sigma}$ and we shall use the
following subgroup
$$
J^{\sigma}=\{\; \; zz^{\sigma}  \quad \vert\quad z \in V({\Bbb
F}_2C), \; zz^{\sigma} \in {\Bbb F}_2C^2\; \;\}.
$$
According to $\sigma$ we distinguish two cases. First let $\sigma$
be the $*$-involution. We shall prove that $J^*=S_*(C)^2$ and its
order is $2^{2^{n-2}-1}$. Indeed, if $u \in J^*$, then $u=zz^*\in
{\Bbb F}_2C^2$ for some $z \in V({\Bbb F}_2C)$ and by Lemma $3$,
$$
u=zz^*=1+\sum_{j=1}^{2^{n-2}-1} \alpha_j(a^{2j}+a^{-2j}).
$$
Now the $*$-symmetric element
$y=1+\sum_{j=1}^{2^{n-2}-1} \alpha_j(a^{j}+a^{-j})$ is such that $y^2=u$.
This shows that  $J^*$ is a subgroup of
$S_*(C)^2$. Conversely, if $u \in S_*(C)^2$, then $u=u^*$ and
there exists $y \in S_*(C)$ with $y^2=u$. Then $yy^*=y^2=u$ and
$u\in J^*$, so $J^*=S_*(C)^2$. The equality of these groups shows
that each $y \in S_*(C)$ satisfies
$$
y^2=yy^*=1+\sum_{j=1}^{2^{n-2}-1} \alpha_j (a^{2j}+a^{-2j}) \in
S_*(C)^2,
$$
hence the order of $S_*(C)^2$ is $2^{2^{n-2}-1}$.

The map  $\psi_*\,:\,V({\Bbb F}_2C) \rightarrow S_*(C)$ induces
the epimorphism $\psi_*\,:\,H_0^* \rightarrow J^*$, and its kernel
is the $*$-unitary subgroup $V_*({\Bbb F}_2C)$ and according to
(2) $V_*({\Bbb F}_2C)$ has order  $2^{2^{n-1}+1}$. Thus
$$
|H_0^*|=|V_*({\Bbb F}_2C)|\cdot|J^*|=|V_*({\Bbb
F}_2C)|\cdot|S_*(C)^2|=2^{3\cdot 2^{n-2}}.
$$

\bigskip
Now, let $\sigma$ be the $\circledast$-involution. We shall show
that
$$
J^{\circledast}=\gp{\; a^{2^{n-1}}}\times S_{\circledast}(C)^2
\quad\text{and}\quad |S_{\circledast}(C)^2|=2^{\frac{|C|}{4}-1}.
$$
For each $y \in S_{\circledast}(C)$ we have
$$
y^2=yy^{\circledast}=1+\sum_{j=1}^{2^{n-2}-1}\beta_{j}(a^{2j}+a^{-2j})
\in S_{\circledast}(C)^2,
$$
and Lemma $4$ asserts that $a^{2^{n-1}} \notin
S_{\circledast}(C)^2$ and $
|S_{\circledast}(C)^2|=2^{\frac{|C|}{4}-1}=2^{2^{n-2}-1}$.

Let  $u \in J^{\circledast}$. Then $u=zz^{\circledast}\in {\Bbb
F}_2C^2$ for some $z \in V({\Bbb F}_2C)$ and again by Lemma $4$ we
obtain that
$$
\split
u=zz^{\circledast}=&\delta_1(1+a^{2^{n-1}})+\\
+&a^{2^{n-1}}\Big(1+
\sum_{i=1}^{2^{n-2}-1}\delta_{j}\big(a^{2j+2^{n-1}}+a^{-(2j+2^{n-1})}\big)
\Big), \quad\quad \delta_i \in {\Bbb F}_2. \endsplit$$

If $\delta_1=0$, then we consider the $\circledast$-symmetric
element
$$
\split
 y_1=1+&\sum_{j\in\{\; 1,3,\dots,2^{n-2}-1\;\}}
\delta_j\big(a^{j+2^{n-2}}+a^{-(j+2^{n-2}+2^{n-1})}\big)\\ +&
\sum_{j\in \{\;
2,4,\dots,2^{n-2}-2\;\}}\delta_j\big(a^{j+2^{n-2}}+a^{-(j+2^{n-2})}\big)
\in S_{\circledast}(C),
\endsplit
$$
which has the property
$$
y^2_1=1+
\sum_{j=1}^{2^{n-2}-1}\delta_{j}\Big(a^{2j+2^{n-1}}+a^{-(2j+2^{n-1})}\Big)
\in S_{\circledast}(C)^2.
$$
Therefore $u=a^{2^{n-1}}y^2_1\in \gp{\; a^{2^{n-1}}}\times
S_{\circledast}(C)^2$.

\noindent Now assume that
$\delta_1=1$. Then the $\circledast$-symmetric element
$$
\split y_2=1+&\sum_{j\in \{\;
1,3,\cdots,2^{n-1}-1\;\}}\delta_j(a^{j}+a^{-j+2^{n-1}})\\+&
\sum_{j\in \{\; 2,4,\cdots,2^{n-1}-2\;\}}\delta_j(a^{j}+a^{-j})
\endsplit
$$
is such that $u=y^2_2$ and $u \in J^{\circledast}\subseteq \gp{\;
a^{2^{n-1}}}\times S_{\circledast}(C)^2$. Conversely, if
\newline $u \in \gp{\; a^{2^{n-1}}}\times S_{\circledast}(C)^2$,
then $u=(a^{2^{n-1}})^td$ for some $d \in S_{\circledast}(C)^2$.
Choose $w \in S_{\circledast}(C)$ with $w^2=d$. Then
$$
\split u=(a^{2^{n-1}})^td=&a^t a^{t(2^{n-1}-1)}w^2=\\=&a^tw
a^{t(2^{n-1}-1)}w=a^tw\cdot(a^tw)^{\circledast} \in
J^{\circledast}. \endsplit
$$

\noindent Consequently, $J^{\circledast}=\gp{\; a^{2^{n-1}}}\times
S_{\circledast}(C)^2$.

Again $\psi_{\circledast}\,:\,V({\Bbb F}_2C) \rightarrow
S_{\circledast}(C)$ induces the epimorphism
$\psi_{\circledast}\,:\,H_0^{\circledast} \rightarrow
J^{\circledast}$ and its kernel is the ${\circledast}$-unitary
subgroup $V_{\circledast}({\Bbb F}_2C)$. This shows that
$$
|H_0^{\circledast}|=|V_{\circledast}({\Bbb
F}_2C)|\cdot|J^{\circledast}|=|V_{\circledast}({\Bbb
F}_2C)|\cdot|\gp{\; a^{2^{n-1}}}\times S_{\circledast}(C)^2|
=2^{3\cdot 2^{n-2}}.
$$

We have established
$$
H_0^{\sigma}\subset H_2^{\sigma}\subset H_4^{\sigma} \subset \dots
\subset H_{2^{n-1}-2}^{\sigma}=V({\Bbb F}_2C),
$$
and $|H_0^{\sigma}|=2^{3\cdot 2^{n-2}}$ from which follows that $
[V({\Bbb F}_2C):H_0^{\sigma}]=2^{2^{n-2}-1}$ and by the second
part of this lemma, the number of different subgroups
$H_{2l}^{\sigma}$ is $2^{n-2}-1$. This is possible if and only if
$[H_{2l+2}^{\sigma}:H_{2l}^{\sigma}]=2$ and we get that the order
of $H_{2l}^{\sigma}$ is $2^{3\cdot 2^{n-2}+l}$ for every
$2l<2^{n-1}$.

\rightline{$\square$}
\enddemo

For  $0\leq i < 2^{n}$  the set
$$
 L_i^{\sigma}=\{\;  h \in V({\Bbb F}_2C)[2] \quad \vert\quad (h+h^{\sigma})(1+a)^{i}=0 \;\}
$$
is a subgroup of
$$
V({\Bbb F}_2C)[2]=\Big\{\; (\sum_{j=0}^{2^{n-1}-1}\alpha_j
a^j)(1+a)^{2^{n-1}}\in V({\Bbb F}_2C) \mid \alpha_j \in  {\Bbb
F}_2\;\Big\}. \tag 12
$$
Indeed,  each $h_k \in L^{\sigma}_{i}$ such that
$h_k(1+a)^i=h_k^{\sigma}(1+a)^i$, so
$$
h_1h_2(1+a)^i=h_1h_2^{\sigma}(1+a)^i=h_1^{\sigma}h_2^{\sigma}(1+a)^i=
(h_1h_2)^{\sigma}(1+a)^i,
$$
and  $L^{\sigma}_{i}$ is a subgroup. Moreover, easy calculations
show that $h^{\circledast}=h^*$  for each  $h \in V({\Bbb
F}_2C)[2]$. Therefore,
$$
L_i^{\circledast}=L_i^*.\tag 13
$$
Consequently it is sufficient to investigate the properties of
$L^*_{i}$.

\proclaim{Lemma 10} The subgroup $L^*_{i}$ has the following properties:
\item{$1.$} $L^*_{i}=V({\Bbb F}_2C)[2]$ for  $i\geq 2^{n-1}$.
\item{$2.$} For even indices $i=2l<2^{n-1}$
the subgroups $L^*_{2l}$ satisfy
$$
L_0^*\subset L_2^* \subset L_4^* \subset\cdots \subset L_{2l}^*
\subset \cdots \subset L_{2^{n-1}-2}^*=V({\Bbb F}_2C)[2],
$$
and the order of $L_{2l}^*$ is $2^{2^{n-2}+1+l}$.
\item{$3.$}  For odd index
$i=2l+1<2^{n-1}$ the subgroup $L_{2l}^*$ coincides with $L_{2l+1}^*$.
\endproclaim

\demo{Proof} First let $i\geq 2^{n-1}-2$. By (12) each $h \in
V({\Bbb F}_2C)[2]$ can be  represented in the form
$$h=
1+\alpha_{2^{n-1}}(1+a)^{2^{n-1}}+\alpha_{2^{n-1}+1}(1+a)^{2^{n-1}+1}+u
$$
for some $u \in A^{2^{n-1}+2}({{\Bbb F}_2}C)$. The formula $(10)$
asserts
$$
h^*=1+\alpha_{2^{n-1}}(1+a)^{2^{n-1}}+\alpha_{2^{n-1}+1}(1+a)^{2^{n-1}+1}+v
$$
for some $v \in A^{2^{n-1}+2}({{\Bbb F}_2}C)$ and $h+h^*=u+v \in
A^{2^{n-1}+2}({{\Bbb F}_2}C)$. But $A^{2^{n-1}+2}({{\Bbb
F}_2}C)\subseteq Ann\big((1+a)^{i}\big)$ for all $i\geq
2^{n-1}-2$, and we conclude that $(h+h^*)(1+a)^{i}=0$. Therefore
$L^*_{i}=V({\Bbb F}_2C)[2]$ for all $i\geq 2^{n-1}$ or
$i=2^{n-1}-2$.
\bigskip

Each symmetric unit $h \in L_0^*$ can
be written as
$$
h=1+\Big(\alpha_0+\alpha_{2^{n-2}}a^{2^{n-2}}+\sum_{j=1}^{2^{n-2}-1}\alpha_j(a^j+a^{-j+2^{n-1}})\Big)
(1+a)^{2^{n-1}},
$$
where $\alpha_j \in {\Bbb F}_2$. This shows that $h$ has
${2^{n-2}+1}$ independent coefficients, thus the order of $L_0^*$
is $2^{2^{n-2}+1}$.

Now we verify that $L_{2l}^*$ is a proper subgroup of
$L_{2(l+1)}^*$ for all $2l<2^{n-1}-2$. According to (10) for
$u=1+(1+a)^{2^n-(2l+3)}$ we have
$$u^*=
1+(1+a)^{2^n-(2l+3)}+(1+a)^{2^n-(2l+2)}+z,
$$
where $z \in A^{2^n-(2l+1)}({{\Bbb F}_2}C)$. This shows that
$u+u^*=(1+a)^{2^n-(2l+2)}+z$, so
$$
(u+u^*)(1+a)^{2l}=(1+a)^{2^n-2}+(1+a)^{2l}z \ne 0,
$$
and $u \notin L_{2l}$. But $(u+u^*)(1+a)^{2l+2}=(1+a)^{2^n}=0$,
this means that $u \in L_{2l+2}^*$ and $L_{2l}^*$ is a proper
subgroup of $L_{2l+2}^*$.

Clearly  $|V({\Bbb F}_2C)[2]|=2^{2^{n-1}}$, and  $ [V({\Bbb
F}_2C)[2]:L_0^*]=2^{2^{n-2}-1}$. As we saw above the number of
different subgroups $L_{2l}^*$ is $2^{n-2}-1$. The only posibility
is that $[L_{2l+2}^*:L_{2l}^*]=2$ and $L_{2l}^*$ has order
$2^{2^{n-2}+1+l}$ for all $2l<2^{n-1}$.

Finally, we use again the element $u=1+(1+a)^{2^n-(2l+3)}$. We
note that
$$
(u+u^*)(1+a)^{2l+1}=(1+a)^{2^n-1} \ne 0,
$$
so $u\not\in
L_{2l+1}^* $ but $u \in L_{2l+2}^*$. It is easy to see that
$L_{2l}^* \subseteq L_{2l+1}^*\subset L_{2l+2}^*$ and
$[L_{2l+2}^*:L_{2l}^*]=2$, so $L_{2l}^*$ coincides with
$L_{2l+1}^*$.

\rightline{$\square$}
\enddemo

\subhead Proof of the main theorem
\endsubhead
\bigskip

We divide the proof into 3 parts. We begin with the dihedral
group. Let us determine the number of the units $x_1+x_2b_2$ of
type $1$ in $V({\Bbb F}_2D_{2^{n+1}})$, where $x_2=0$ or
$x_2=\gamma (1+a)^i$, $\gamma$ is a unit and $i>0$.

If $x_2=0$, then  by $(5)$ we have  $x_1^2=1$. Therefore the
number of units of type $1$ with $x_2=0$ coincides with the order
of $V({\Bbb F}_2C)[2]$.

Now let $0<i<2^{n-1}$ and $x_2=\gamma(1+a)^i$. Then
$$
x_1^2=\gamma\gamma^*(a+a^{-1})^i+1; \tag 14
$$
$$
(x_1+x_1^*)(1+a)^i=0. \tag 15
$$
By $(14)$ the element $\gamma\gamma^*(a+a^{-1})^i$ belongs to
${\Bbb F}_2C^2$, further according to Lemma $8$ the number $i=2l$
is even and
$$
H_{2l}^*=\{\; \gamma\in V({\Bbb F}_2C)\quad\vert\quad
\gamma\gamma^*(a+a^{-1})^{2l} \in  {\Bbb F}_2C^2\;\}
$$
is a subgroup of $V({\Bbb F}_2C)$.

For fixed $\gamma$  and $i$ we determine the number of units
$x_1+x_2b_2$ of type $1$, which satisfy the  conditions $(14)$ and
$(15)$. If the unit $x_1'+x_2b$ is also type $1$, then
$x_1'x_1^{-1}$ is a unit of order two and
$$
(x_1'x_1^{-1}+(x_1'x_1^{-1})^*)(1+a)^{2l}=
x_1'(x_1^{-1}+(x_1^{-1})^*)(1+a)^{2l}=0.
$$
Therefore
$$
x_1'x_1^{-1}\in L_{2l}^*=\{\; h\in V({\Bbb
F}_2C)[2]\quad\vert\quad (h+h^*)(1+a)^{2l}=0\;\},
$$
so $x_1' \in x_1\cdot L_{2l}^*$ and the number of different
elements $x_1'$ is $|L_{2l}^*|$.

Finally,  we shall determine the number of elements
$x_2=\gamma(1+a)^i$ for a fixed $i$, where $\gamma \in V({\Bbb
F}_2C)$. This number coincides with  the cardinality of the set
$$
K_{2l}=\{\; \gamma(1+a)^{2l}\;|\; \gamma \in V({\Bbb F}_2C),\;
\gamma\gamma^*(a+a^{-1})^{2l}\in {\Bbb F}_2C^2 \;\}.
$$
Clearly $K_{2l}$ coincides with $H_i^*(1+a)^{2l}$. But
$\gamma(1+a)^{2l}=\gamma'(1+a)^{2l}$ if and only if
$\gamma^{-1}\gamma'(1+a)^{2l}=(1+a)^{2l}$, so
$\gamma^{-1}\gamma'\in S_{2l}$. We have established the equality
$|K_{2l}|=\frac{|H_{2l}^*|}{|S_{2l}|}$. Hence the number of units
of type $1$ with the form $x_1+x_2b_2$ for a fixed $0<2l <2^{n-1}$
is equal to
$$
\textstyle \frac{|H_{2l}^*|}{|S_{2l}|}\cdot|L_{2l}^*|= 2^{2^n+1}.
$$

Now consider the case, when    $2^{n-1}\leq i< 2^n$. Then
$(a+a^{-1})^i=0$, and (14) implies that the unit $x_1+x_2b_2$ is
such that  $x_1\in V({\Bbb F}_2C)[2]$ and $(15)$ is always
satisfied. Thus $x_2\in V({\Bbb F}_2C)(1+a)^i$ and there are
$$
\textstyle \frac{|V({\Bbb F}_2C)|}{|S_i|}\cdot |V({\Bbb F}_2C)[2]|
$$
different  units of type $1$ in $V({\Bbb F}_2D_{2^{n+1}})$. We
obtain that in $V({\Bbb F}_2D_{2^{n+1}})$ the number of units of
type $1$ is  equal to
$$
\textstyle |V({\Bbb F}_2C)[2]|+\sum_{l=1}^{2^{n-2}-1}
\frac{|H_{2l}^*|}{|S_{2l}|}\cdot|L_{2l}^*|
+\sum_{j=2^{n-1}}^{2^n-1} \frac{|V({\Bbb F}_2C)|}{|S_j|}\cdot
|V({\Bbb F}_2C)[2]|=
$$
$$
\textstyle 2^{2^{n-1}}+ \sum_{j=1}^{2^{n-2}-1} 2^{2^n+1} +
2^{2^n-1}\sum_{j=2^{n-1}}^{2^{n}-1} 2^{2^{n-1}}2^{-j}=
2^{{2^n}}(2^{n-1}-1).
$$

Now let us consider the number  of units of type $2$ in $V({\Bbb
F}_2D_{2^{n+1}})$. If $x_1+x_2b_2$ is a unit of type $2$, then
$x_2$ is a unit and according to $(5)$ we have  $x_1=x_1^*$ and
$x_2x_2^*=(1+x_1)^2$. Evidently $1+x_1$ is a $*$-symmetric unit
and $x_2x_2^*\in V({\Bbb F}_2C^2)$. For a fixed unit $x_2$, by the
first part of Lemma $2$ the set
$$
\{\;  1+x_1 \in S_*(C)\quad\vert\quad (1+x_1)^2=x_2x_2^* \;\}
$$
is a coset of $S_*(C)$ by $S_*(C)[2]$. Therefore the number of the
different $x_1$ is $|S_*(C)[2]|$. Clearly $L_0^*$ coincides with
$S_*(C)[2]$, so for a fixed $x_2$ the number of the different
$x_1$ is $|L_0^*|$.

Since $H_0^*=\{\; h \in V({\Bbb F}_2C) \; | \; hh^*\in V({\Bbb
F}_2C^2) \;\}$, the number of the different $x_2$ coincides with
$|H_0^*|$. This shows that in $V({\Bbb F}_2D_{2^{n+1}})$  the
number of the units of type $2$ is
$$
|H_0^*|\cdot|L_0^*|=2^{3\cdot 2^{n-2}+2^{n-2}+1}=2^{2^{n}+1}
$$
and
$$
\Theta_{D_{2^{n+1}}}(2)=2^{2^n+n-1}+2^{2^n}.
$$

Now let us consider the generalized quaternion group. Let
$x_1+x_2b_1$ be  a  unit of type $2$. Similarly to the previous
case,  $x_1$ is  $*$-symmetric  and
$$
x_1^2+x_2x_2^*a^{2^{n-1}}=x_1x_1^*+x_2x_2^*a^{2^{n-1}}=1.
$$
Since $1 \not\in supp\,(x_1x_1^*)$ and $1 \not\in
supp\,(x_2x_2^*a^{2^{n-1}})$ by Lemma $3$ this equality is
impossible. Therefore, there is no unit of type $2$ in $V({\Bbb
F}_2Q_{2^{n+1}})$.

By Lemma $7$ in $V({\Bbb F}_2Q_{2^{n+1}})$ the number of units of
type $1$ is $(2^{n-1}-1)2^{2^n}$ and
$$
\Theta_{Q_{2^{n+1}}}(2)=2^{2^n+n-1}-2^{2^n}.
$$

\bigskip
Finally, let us consider the semidihedral group. According to
$(13)$ and Lemma $8$, we have $|H_i^{\circledast}|=|H_i^*|$ and
$|L_i^{\circledast}|=|L_i^*|$ for each $i<2^n$. Thus in $V({\Bbb
F}_2SD_{2^{n+1}})$ the number of units of type $1$ is
$(2^{n-1}-1)2^{2^n}$ as in the group of units $V({\Bbb
F}_2D_{2^{n+1}})$ for the dihedral group.

Now  we shall determine the number of  units $x_1+x_2b_3  \in
V({\Bbb F}_2SD_{2^{n+1}})$ of type $2$. Then $x_2$ is a unit and
from $(5)$ it follows that $x_2^{\circledast}x_2=(x_1+1)^2\in
V({\Bbb F}_2C^2)$ and $(x_1+x_1^{\circledast})x_2=0$. Since $x_2$
is a unit we get that $x_1$ is a $\circledast$-symmetric,  so
$1+x_1$ is a $\circledast$-symmetric unit as well. Evidently
$x_2^{\circledast}x_2=(1+x_1)^2 \in S_{\circledast}(C)^2$. We have
seen that
$$
J_{\circledast}=\{\;  x_2x_2^{\circledast}  \; |\; x_2 \in V({\Bbb
F}_2C), \; x_2x_2^{\circledast} \in {\Bbb F}_2C^2 \;\}=\gp{\;
a^{2^{n-1}}} \times S_{\circledast}(C)^2,
$$
thus the number of different $x_2$ is
$$
\textstyle \Big|\{\;\; x_2 \in V({\Bbb F}C) \;\vert\;
x_2x_2^{\circledast} \in S_{\circledast}(C)^2\;\}\Big|=
$$$$
\textstyle =|V_{\circledast}({\Bbb F}C)|\cdot
|S_{\circledast}(C)^2|=|V_{\circledast}({\Bbb F}C)|\cdot
\frac{|J_{\circledast}|}{2}=\frac{|H_0^{\circledast}|}{2}.
$$

 For a fixed unit $x_2$ the set $\{\;  1+x_1 \in S_{\circledast}(C) \quad\vert\quad
(1+x_1)^2=x_2x_2^{\circledast} \;\}$ is a coset of
$S_{\circledast}(C)$ by $S_{\circledast}(C)[2]$ as the first part
of Lemma $2$ asserts. Clearly $L_0^{\circledast}$ coincides with
$S_{\circledast}(C)[2]$ so for a fixed $x_2$ the number of the
different $x_1$ is $|L_0^{\circledast}|$. We obtain that the
number of units of type $2$ is
$$
\textstyle \frac{|H_0^{\circledast}|}{2}\cdot
|L_0^{\circledast}|=2^{3\cdot 2^{n-2}+2^{n-2}}=2^{2^{n}}
$$
and $ \Theta_{SD_{2^{n+1}}}(2)=2^{2^n+n-1}$, the proof is done.

\rightline{$\square$}

\enddocument